%dmp.tex:Using GENERATINGFUNCTIONOLOGY to Enumerate Distinct-Multiplicity Partitions
%%a Plain TeX file by Doron Zeilberger (x pages)

%begin macros

\baselineskip=14pt
\parskip=10pt

\font\eightrm=cmr8 
\font\eighttt=cmtt8
\magnification=\magstephalf
\def\F{{\cal F}}
\def\P{{\cal P}}

\def\1{{\overline{1}}}
\def\2{{\overline{2}}}
\parindent=0pt
\overfullrule=0in

\def\frac#1#2{{#1 \over #2}}
%\headline={\rm  \ifodd\pageno  \RightHead  \else  \LeftHead  \fi}
%\def\RightHead{\centerline{
%Title
%}}
%\def\LeftHead{ \centerline{Doron Zeilberger}}
%end macros
\bf
\centerline
{
Using GENERATINGFUNCTIONOLOGY to Enumerate Distinct-Multiplicity Partitions
}
\rm
\bigskip
\centerline{ {\it
Doron 
ZEILBERGER}\footnote{$^1$}
{\eightrm  \raggedright
Department of Mathematics, Rutgers University (New Brunswick),
Hill Center-Busch Campus, 110 Frelinghuysen Rd., Piscataway,
NJ 08854-8019, USA.
%\break
{\eighttt zeilberg  at math dot rutgers dot edu} ,
\hfill \break
{\eighttt http://www.math.rutgers.edu/\~{}zeilberg/} .
Jan. 18, 2012.
Accompanied by Maple package {\eighttt DMP} 
downloadable from {\eighttt http://www.math.rutgers.edu/\~{}zeilberg/tokhniot/DMP} .
Sample input and output files may be viewed in the front of this article:
{\eighttt http://www.math.rutgers.edu/\~{}zeilberg/mamarim/mamarimhtml/dmp.html}  .
Supported in part by the USA National Science Foundation.
}
}

{\it In fond memory of Guru Herbert Saul WILF (28 Sivan 5691- 12 Tevet 5772) zecher gaon l'bracha}

{\bf Preamble}

About a year ago, Herb Wilf[W1] posed, on-line, eight intriguing problems.
I don't know the {\it answer} to any of them, but I will say something about
the sixth question.

{\bf Herb Wilf 6th Question}:
{\it Let $T(n)$ be the set of partitions of $n$ for which the (nonzero) multiplicities of its parts are
all different, and write $f(n)=|T(n)|$. See Sloane's sequence {\bf A098859} for a table of values.
Find any interesting theorems about $f(n)$ \dots }

First, I will explain how to compute the first few terms of $f(n)$. Shalosh can easily get
the first $250$ terms, but as $n$ gets larger it gets harder and harder to compute,  unlike its unrestricted cousin, $p(n)$.
I conjecture that the fastest algorithm takes exponential time, but I have no idea how 
to prove that claim. I am impressed that, according to Sloane, Maciej Ireneusz Wilczynsk computed $508$ terms.

Recall that the generating function for the number of integer partitions of $n$ whose largest part is $\leq m$, $p_m(n)$, is
the very simple rational function
$$
\sum_{n=0}^{\infty} p_m(n)q^n \, = \, \frac{1}{(1-q) (1-q^2) \cdots (1-q^m) } \quad  .
$$
The main purpose of this note is to describe, using {\it Generatingfunctionology}, so vividly and lucidly
{\it preached} in W's {\bf classic book} [W2],
how to compute the generating function (that also turns out to be rational)
for the number of partitions of $n$ whose largest part is $\leq m$ {\bf and} all its (nonzero) multiplicities are distinct,
let's call it $f_m(n)$.
As $m$ gets larger, the formulas get more and more complicated, but we sure do have an {\bf answer},
in the sense of the {\bf classic article} [W3], for any {\it fixed} $m$, but of course not for a {\it symbolic} $m$.

Even more is true! Because, like $\frac{1}{(1-q) (1-q^2) \cdots (1-q^m) }$, the generating function
of $f_m(n)$, $\sum_{n=0}^{\infty} f_m(n)q^n$, turns out (as we will see) 
to only have roots-of-unity poles, whose highest order is $m$,
it follows that $f_m(n)$ is a {\it quasi-polynomial} of degree $m-1$ in $n$.
Now that's a very good answer! (in W's sense, albeit only for a fixed $m$).

\vfill
\eject

{\bf How to Compute Many terms of $f(n)$}?

$p_m(n)$ is very easy to compute. For example, one may use  the recurrence
$$
p_m(n)=p_{m-1}(n)+\sum_{i=1}^{\lfloor n/m \rfloor} p_{m-1}(n-mi) \quad ,
$$
together with the {\it initial condition} $p_1(n)=1$, $p_m(0)=1$.

How can we adapt this in order to compute $f_m(n)$? 
The contribution from the partitions counted by $f_m(n)$ where  $m$ does not show up is $f_{m-1}(n)$,
in analogy with the $p_{m-1}(n)$ term in the above recurrence. But if $m$ {\it does} show up, it does so
with a certain multiplicity, $i$, say, where $1 \leq i \leq \lfloor n/m \rfloor $, and
removing these $i$ copies  of $m$ results in a partition counted by $f_{m-1}(n-mi)$
-so all its multiplicities are different- and {\bf in addition} none of these multiplicities may be $i$.
Continuing, we are forced to introduce a much more general discrete function $f_m(n;S)$
whose arguments are $m$ and $n$ and a set of ``forbidden multiplicities'', $S$.

So let's define $f_m(n;S)$ to be the number of  partitions of $n$ with parts $\leq m$, with all its
multiplicities distinct {\bf and} none of these multiplicities belonging to $S$.
Our intermediate object of desire, $f_m(n)$, is simply $f_m(n; \emptyset)$, and the ultimate object,
$f(n)$, is $f_n(n; \emptyset)$.

The recurrence for $f_m(n;S)$ is, naturally

$$
f_m(n; S)=f_{m-1}(n; S)+\sum_{i=1, i \not \in S}^{\lfloor n/m \rfloor} f_{m-1}(n-im; S \cup \{i\}) \quad ,
$$
because once we decided on the number of times $m$ shows up, let's call it $i$,
where $i$ is between $1$ and $\lfloor m/m \rfloor$ and $i \not \in S$,
the partition (of $n-mi$) obtained by removing these $i$ copies of $m$ must forbid
the set of multiplicities $S \cup \{ i \}$.

In the Maple package {\tt DMP}, procedure {\tt qnmS(n,m,S)} implements $f_m(n;S)$ and
procedure {\tt qn(n)} implements $f(n)$.

{\bf Inclusion-Exclusion}

Let $\P_m(n)$ be the set of partitions of $n$ whose parts are all $\leq m$,
in other words, the set that $p_m(n)$ is counting.
Consider the set of {\it all} partitions whose largest part is $\leq m$,
where we write a partition in {\it frequency notation}:
$$
\P_m := \{ 1^{a_1} 2^{a_2} \dots  m^{a_m} \, | \, a_1 ,\dots, a_m \geq 0 \} \quad .
$$
For example $1^3 2^5 4^2$ is the partition of twenty-one usually written as $4422222111$.
Introducing {\it symbols} $x_1,x_2, \dots ,  x_m$, we define the
{\it Weight} of a partition to be
$$
Weight(1^{a_1} 2^{a_2} \dots  m^{a_m}):=x_1^{a_1} x_2^{a_2} \cdots x_m^{a_m} \quad .
$$
The {\it weight-enumerator} of $\P_m$ is, by {\bf ordinary}-generatingfunctionology
$$
Weight(\P_m)=\frac{1}{(1-x_1)(1-x_2) \cdots (1-x_m)} \quad ,
$$
since we make $m$ {\it independent} decisions: 

$\bullet$ how many copies of $1$?(Weight enumerator  $=1+x_1+x_1^2 + \dots =(1-x_1)^{-1}$) ,

$\bullet$ how many copies of $2$? (Weight enumerator $=1+x_2+x_2^2 + \dots =(1-x_2)^{-1}$) ,

$\dots$

$\bullet$ how many copies of $m$? (Weight enumerator $=1+x_m+x_m^2 + \dots =(1-x_m)^{-1}$).

But we want to find the weight-enumerator of the much-harder-to-weight-count set

$$
\F_m := \{ 1^{a_1} 2^{a_2} \dots  m^{a_m} \, | \, a_1 ,\dots, a_m \geq 0 \,  ;  \, a_i \neq a_j \,\,\, (if \quad a_i>0, a_j > 0) \} \quad .
$$

Calling the members of $\F_m$ {\it good}, we see that a member of $\P_m$ is good if it does
{\bf not} belong to any of the following ${{m} \choose {2}}$ sets, $S_{ij}$ $1 \leq i < j \leq m$:
$$
S_{ij}:= \{ 1^{a_1} 2^{a_2} \dots  m^{a_m} \in \P_m \, | a_i=a_j>0\} \quad .
$$
By inclusion-exclusion, the weight-enumerator of $\F_m$ is 
$$
\sum_{G} (-1)^{|G|} Weight \left ( { { \cap} \atop {ij \in G}} \,\, S_{ij} \right ) \quad ,
$$
where the summation ranges over {\bf all} $2^{m(m-1)/2}$ subsets of $\{ (i,j) \, \vert \, 1 \leq i <j \leq m \}$.

But the $G$'s can be naturally viewed as {\it labeled  graphs} on $m$ vertices. Such a graph has
several connected components, and together they naturally induce a {\it set partition} 
$\{ C_1, C_2, \dots , C_r \}$ of $\{1, 2, \dots, m\}$. 
We have:
$$
Weight \left ( { { \cap} \atop {ij \in G}} \,\, S_{ij} \right )
\, = \, \prod_{i=1}^{r}  weight(C_i) \quad ,
$$
where if $|S|=1$,  $S=\{s\}$, say,  then $weight(S)=\frac{1}{1-x_s}$, and if $|S|=d>1$, 
$S=\{s_1, s_2, \dots s_d\}$, say, then
$$
weight(S)=\frac{x_{s_1} x_{s_2} \cdots x_{s_d}}{1-x_{s_1} x_{s_2} \cdots x_{s_d}} \quad .
$$
To justify the latter, note that if vertices $s_1, s_2, \dots s_d$ all belong to the same connected
component of our graph then, by {\it transitivity}, we have
that all $a_{s_1}=a_{s_2}= \dots =\dots a_{s_d}>0$, and the weight-enumerator is the infinite geometric
series
$$
\sum_{\alpha=1}^{\infty} (x_{s_1} \cdots x_{s_d})^{\alpha}=\frac{x_{s_1} x_{s_2} \cdots x_{s_d}}{1-x_{s_1} x_{s_2} \cdots x_{s_d}} \quad .
$$

But quite a few graphs correspond to any one set-partition. To find out the coefficients in front,
for any set-partition $\{ C_1, C_2, \dots , C_r \}$  of $\{1, \dots , m\}$ we must find
$$
\sum_G (-1)^{|G|} \quad ,
$$
summed over all the graphs that gives rise to the above set partition.
But this is the product of the analogous sums where one focuses on one connected component at a time,
and then multiplies everything together.

Let's digress and figure out $\sum_G (-1)^{|G|}$ over all {\it connected  labeled graphs} on $n$ vertices.
For the sake of clarity, let's, more generally, figure out $\sum_G y^{|G|}$ with a general variable $y$.

By {\bf exponential}-generatingfunctionology[W2] (see also [Z]), this sum is nothing but the coefficient of $t^n/n!$ in
$$
\log \left ( \sum_{i=0}^{\infty} (1+y)^{{{i} \choose {2}}} \frac{t^i}{i!} \right ) \quad .
$$
Going back to $y=-1$, we see that we need the coefficient of $t^n/n!$ in
$$
\log \left ( \sum_{i=0}^{\infty} (1+(-1))^{{{i} \choose {2}}} \frac{t^i}{i!} \right ) =
\log \left ( \sum_{i=0}^{\infty} 0^{{{i} \choose {2}}} \frac{t^i}{i!} \right ) =
\log  (1+t) 
$$
$$
=\sum_{n=1}^{\infty} (-1)^{n-1} \frac{t^n}{n}=
\sum_{n=1}^{\infty} (-1)^{n-1} (n-1)! \frac{t^n}{n!} \quad .
$$
So the desired sum is $(-1)^{n-1} (n-1)!$.

Let's define for any set of positive integers, $S$,
$$
mishkal(S)=
\cases{
1/(1-x_s) ,& if $|S|=1$ where $S=\{s\}$ ;\cr
(-1)^{d-1}(d-1)!(x_{s_1} \cdots x_{s_d})/(1-x_{s_1} \cdots x_{s_d}) ,&  if $|S|=d>1$ where $S=\{s_1, \dots s_d \}$.\cr}   
$$

For any set partition $C=\{C_1, \dots , C_r\}$ let's define
$$
Mishkal(C)= mishkal(C_1) \cdots mishkal(C_r) \quad .
$$
It follows that the weight-enumerator of $\F_m$ according to  $Weight(1^{a_1} 2^{a_2} \dots  m^{a_m}):=x_1^{a_1} x_2^{a_2} \cdots x_m^{a_m}$
is 
$$
\sum_C Mishkal(C) \quad ,
$$

where the sum has $B_m$ terms ($B_m$ being the Bell numbers), one for each set-partition of $\{ 1, \dots , m\}$.

{\bf Finally}, to get an ``explicit'' formula (as a sum of $B_m$ terms, each a simple rational function of $q$),
for the generating function $\sum_{n=0}^{\infty} f_m(n)q^n$, all we need is replace $x_i$ by $q^i$, for $i=1 \dots m$, 
getting
$$
\sum_{n=0}^{\infty} f_m(n)q^n= \sum_C Poids(C) \quad ,
$$
where for a set partition $C=\{C_1, \dots , C_r\}$ 
$$
Poids(C)= poids(C_1) \cdots poids(C_r) \quad ,
$$
and where for an individual set $S$:
$$
poids(S)=
\cases{
1/(1-q^s) ,& if $|S|=1$ where $S=\{s\}$ ;\cr
(-1)^{d-1} (d-1)!q^{s_1+ \dots + s_d}/(1-q^{s_1+ \dots + s_d}) ,&  if $|S|=d>1$ where $S=\{s_1, \dots s_d \}$.\cr}   \quad  .
$$

It follows that indeed $f_m(n)$ is a quasi-polynomial of degree $m-1$ in $n$.
Furthermore, since the only pole that has multiplicity $m$ is $q=1$, it follows that
the leading term (of degree $m-1$) is a pure polynomial.

The generating function, $\sum_{n=0}^{\infty} f_m(n)q^n$, for any desired positive integer $m$, 
is implemented in procedure {\tt GFmq(m,q)} in the Maple package {\tt DMP}.
For the Weight-enumerator (or rather with $x_i$ replaced by $q^i x_i$, for $i=1, \dots ,m$), 
see {\tt GFmxq(m,x,q)}. Since the Bell numbers grow very fast, the formulas get complicated rather fast, but
{\it in principle} we do have a very nice {\it answer} for any specific $m$, but {\it in practice}, for large
$m$ it is only ``nice'' in principle. Of course it is {\it anything but} nice when viewed also as function of
$m$, and that's why $f(n)=f_n(n)$ is {\bf probably} very hard  to compute for larger $n$.

To see the outputs of {\tt GFmq(m,q)} for $1 \leq {\tt m} \leq 8$ see: \hfill\break
{\tt http://www.math.rutgers.edu/\~{}zeilberg/tokhniot/oDMP3 } .

{\bf Asymptotics}

 Recall that  Hardy and Ramanujan tell us
that as $n$ goes to infinity, $p(n)$ is asymptotic 
to $\frac{1}{4n\sqrt{3}}exp(C\sqrt{n})$ where $C=\sqrt{2/3} \pi=2.565099661\dots$,
and hence $\log p(n)/\sqrt{n}$ converges to $C$. By looking at the
sequence $\log f(n)/\sqrt{n}$ for $1 \leq n \leq 508$,  it seems that this
too converges to a limit, that appears to be 
a bit larger than $1.517$ (but of course way less than $2.565099661\dots$).
Let's call that constant the {\it Wilf constant}.

The numerical evidence is here:  {\tt http://www.math.rutgers.edu/\~{}zeilberg/tokhniot/oDMP4 } .

Let me conclude with two challenges.

$\bullet$ Prove that the Wilf constant exists.

$\bullet$ Determine the exact value of the Wilf constant (if it exists) in terms of $\pi$ or other famous constants. 
Failing this, find non-trivial rigorous lower and upper bounds.

\vfill\eject

{\bf References}

[W1] Herbert Wilf, {\it Some Unsolved problems}, \hfill\break
{\tt http://www.math.upenn.edu/\%7Ewilf/website/UnsolvedProblems.pdf}, posted: Dec. 13, 2010.
(viewed Jan. 16, 2012)

[W2] Herbert S. Wilf, {\bf ``Generatingfunctionology''}, Academic Press, First edition 1990, Second Edition 1994.
Third Edition: AK Peters, 2005. Second edition is freely downloadable from \hfill\break
{\tt http://www.math.upenn.edu/\%7Ewilf/gfology2.pdf} .

[W3] Herbert S. Wilf, {\it What is an Answer?}, Amer. Math. Monthly {\bf 89} (1982), 289-292.

[Z] Doron Zeilberger, {\it Enumerative and Algebraic Combinatorics}, in:
{\bf ``Princeton Companion to Mathematics''} , (Timothy Gowers, ed.), Princeton University Press, 550-561.
Available from: \hfill\break
{\tt http://www.math.rutgers.edu/\~{}zeilberg/mamarim/mamarimPDF/enu.pdf}
\end